\documentstyle[12pt]{article}
\newcommand{\ptl}{\partial}

\newcommand{\G}{\Gamma}

\newcommand{\vs}{\varsigma}

\newcommand{\ol}{\overline}

\def\DD{{\cal D}}
\def\AA{{\cal A}}
\def\WW{{\cal W}}

\def \Z{\hbox{$Z\hskip -5.2pt Z$}}

\def\sZ{\hbox{$\sc Z\hskip -4.2pt Z$}}

\def \C{\hbox{$C\hskip -5pt \vrule height 6pt depth 0pt \hskip 6pt$}}

\def\qed{\hfill \hfill \ifhmode\unskip\nobreak\fi\ifmmode\ifinner
         \else\hskip5pt\fi\fi
 \hbox{\hskip5pt\vrule width4pt height6pt depth1.5pt\hskip 1 pt}}
\def\a{\alpha}
\def\b{\beta}
\def\d{\delta}

\def\g{\gamma}
\def\G{\Gamma}
\def\l{\lambda}

\def\sc{\scriptstyle}
\def\ssc{\scriptscriptstyle}
\def\F{\hbox{$I\hskip -4pt F$}}\def \Z{\hbox{$Z\hskip -5.2pt Z$}}
\def\dis{\displaystyle}
\def\cl{\centerline}

\def\nl{\newline}
\def\ol{\overline}

\def\nob#1{$\mbox{#1}$}
\def\wh{\widehat}
\def\rar{\rightarrow}

\def\Rla{\Leftrightarrow}

\def\bs{\backslash}

\def\rb{\raisebox}

\def\vs{\vspace*}
\def\vsp{\vs{-5pt}}
\def\vms{\vs}

\def\ni{\noindent}
\def\hi{\hangindent}
\def\ha{\hangafter}

\begin{document}
\ni
2-COCYCLES ON THE LIE ALGEBRAS OF GENERALIZED\nl \mbox{DIFFERENTIAL} OPERATORS
\par\
\par\ni
                 Yucai Su
\par\ni
                 Department of Applied Mathematics,
                 Shanghai Jiaotong University,
                 Shanghai 200030, China
                 (Email: yucai\_su@cicma.concordia.ca)
\par\
\par\ni
                 ABSTRACT.
In a recent paper by Zhao and the author, the Lie algebras
$\AA[\DD]=\AA\otimes\F[\DD]$ of Weyl type were defined and studied, where
$\AA$ is a commutative associative algebra with an identity element over a
field $\F$ of any characteristic, and $\F[\DD]$ is the polynomial algebra of
a commutative derivation subalgebra $\DD$ of $\AA$. In the present paper,
the 2-cocycles of a class of the above Lie algebras $\AA[\DD]$ (which are
called the Lie algebras of generalized differential operators in the present
paper), with $\F$ being a field of characteristic $0$, are determined. Among
all the 2-cocycles, there is a special one which seems interesting. Using
this 2-cocycle, the central extension of the Lie algebra is defined.
\par\
\par
\cl{             \S1. Introduction
}
\par
We start with a brief definition.
For a Lie algebra $L$ over a field $\F$ of characteristic zero,
a {\it 2-cocycle} on $L$ is a $\F$-bilinear function $\psi:L\times L\rar
\F$ satisfying the following conditions:
$$
\matrix{
\psi(v_1,v_2)=-\psi(v_2,v_1),
\vs{4pt}\hfill\cr
\psi([v_1,v_2],v_3)+\psi([v_2,v_3],v_1)+\psi([v_3,v_1],v_2)=0,
\hfill\cr}
\eqno\matrix{(1.1)\vs{4pt}\cr(1.2)\cr}$$
for $v_1,v_2,v_3\in L$.
Denote by $C^2(L,\F)$ the vector space of 2-cocycles on $L$.
For any $\F$-linear function $f:L\rar\F$, one can define a 2-cocycle
$\psi_f$ as follows
$$
\psi_f(v_1,v_2)=f([v_1,v_2]),
\eqno(1.3)$$
for $v_1,v_2\in L.$
Such a 2-cocycle is called a {\it 2-coboundary} or a {\it trivial
2-cocycle} on $L$.
Denote by $B^2(L,\F)$ the
vector space of 2-coboundaries on $L$.
A 2-cocycle $\phi$ is said to be
{\it equivalent to} a 2-cocycle $\psi$ if $\phi-\psi$ is trivial.
The quotient space
$$
\matrix{
H^2(L,\F)\!\!\!\!&=C^2(L,\F)/B^2(L,\F)
\vs{4pt}\hfill\cr&=
\{\mbox{the equivalent classes of 2-cocycles}\},
\hfill\cr}
\eqno(1.4)$$
is called the {\it 2-cohomology group} of $L$.
\par
The 2-cocycles on Lie algebras play important roles in the central extensions
of Lie \nob{algebras.}
It is well known that all 1-dimensional central extensions of $L$ are
determined by the \nob{2-cohomology} group of $L$.
Central extensions are often used in the structure
\nob{theory}
and the \nob{representation} theory of Kac-Moody algebras [4,5].
Using central extension, we can \nob{construct} many infinite dimensional Lie
algebras, such as affine Lie algebras, infinite \nob{dimensional} \nob{Heisenberg}
algebras, and generalized Virasoro and super-Virasoro algebras,
which have a profound \nob{mathematical} and physical background
[2-4,6,7,12,15-17].
We can describe the structures and some of the representations of
these Lie algebras by using central extension method [5,8].
Since the \nob{cohomology} groups are closely related to the
structures of Lie \nob{algebras,}
the computation of \nob{cohomology} groups
seems to be important and interesting as well.
Berman [1] and Su [13] gave some computation of
the low dimensional cohomology groups for some infinite dimensional Lie
algebras.
\par
Below, we shall introduce the Lie algebras we are going to consider in this
paper.
For any positive integer $n$, an additive subgroup $G$ of $\F^n$ is called
{\it nondegenerate} if $G$ contains an \F-basis of $\F^n$.
Let $\ell_1,...,\ell_4$ be four nonnegative integers such that
$\ell=\sum_{i=1}^4\ell_i>0$. For convenience, we
\vsp
denote
$$
\ell'_i=\sum_{p=1}^i\ell_p,\ \
\ell''_i=\sum_{p=i}^4\ell_p,\ \ i=1,2,3,4.
\vsp
\eqno(1.5)$$
For any $m,n\in\Z$, we use the notation
$$
\ol{m,n}=\{m,m+1,...,n\}\mbox{ if }m\le n.
\eqno(1.6)$$
Take an additive subgroup $\G$ of $\F^{\ell}$
such that for any $\a\in\G$,
$$
\a=(\a_1,\a_2,...,\a_\ell)=(0,...,0,\a_{\ell_1+1},...,\a_\ell),
\eqno(1.7)$$
and such that
$\G$ is nondegenerate as a subgroup of $\F^{\ell''_2}$. Set
$$
\vec J=\Z_+^{\ell'_2}\times\Z^{\ell_3}\times\{0\}^{\ell_4}
\subset\Z^\ell.
\eqno(1.8)$$
Elements in $\Z^\ell$ will be written as
$$
\vec i=(i_1,i_2,...,i_\ell),
\mbox{ and } \vec i=(i_1,...,i_{\ell'_3},0,...,0)
\mbox{ if }\vec i\in\vec J.
\eqno(1.9)$$
For any $i\in\Z$, $p\in\ol{1,\ell}$, we
\vms{-8pt}
denote
$$
i_{[p]}=(0,...,\rb{8pt}{\mbox{$^{^{\sc p}}_{_{\dis i}}$}},0,...,0)\in\Z^\ell.
\vms{-3pt}
\eqno(1.10)$$
Let $\AA$ be the semi-group algebra
$\F[\G\times\vec J]$ with basis $\{x^{\a,\vec i}\,|\,(\a,\vec i)
\in\G\times\vec J\}$ and the algebraic operation
$$
x^{\a,\vec i}\cdot x^{\b,\vec j}=x^{\a+\b,\vec i+\vec j},
\eqno(1.11)$$
for $(\a,\vec i),(\b,\vec j)\in\G\times\vec J.$
Denote the identity element $x^{0,\vec0}$ by 1, and for convenience,
denote
$x^{\a}=x^{\a,\vec 0},\,
t^{\vec i}=x^{\vec 0,\vec i}$ for $\a\in\G,\,\vec i\in\vec J$.
Define the linear transformations
$\{\partial^-_1,...,\partial^-_{\ell'_3},\partial^+_{\ell_1+1},
...,\partial^+_{\ell}\}$ on $A$ by
$$
\partial^-_p(x^{\a,\vec i})=i_px^{\a,\vec i-1_{[p]}},
\ \ \partial^+_q(x^{\a,\vec i})=\a_q x^{\a,\vec i},
\eqno(1.12)$$
for all $p\in\ol{1,\ell'_3},\,q\in\ol{\ell_1+1,\ell}$.
The operators $\partial_p^-$ are called {\it down-grading operators} and
$\partial^+_q$ are called {\it grading operators}.
Set
$$
\partial_p=\partial^-_p,\,\partial_q=\partial^-_q+\partial^+_q,\,
\partial_r=\partial^+_r,
\eqno(1.13)$$
for $p\in\ol{1,\ell_1},\,
q\in\ol{\ell_1\!+\!1,\ell'_3},\,
r\in\ol{\ell'_3+1,\ell}.$
Note that $\partial_p$ is locally nilpotent if $p\in\ol{1,\ell_1}$;
locally finite if $p\in\ol{\ell_1+1,\ell'_2}$; not locally finite if
$p\in\ol{\ell'_2+1,\ell'_3}$; semi-simple if $p\in\ol{\ell'_3+1,\ell}$.
\par
Denote $\DD=\sum_{p=1}^\ell\F\ptl_p$.
The \F-vector
\vsp
space
$$
W=W(\ell_1,\ell_2,\ell_3,\ell_4,\G)=\AA\otimes\DD=\sum_{p=1}^\ell \AA
\partial_p,
\vms{-10pt}
\eqno(1.14)$$
forms a Lie algebra under the bracket
$$
[u\partial_p,v\partial_q]=u\partial_p(v)\partial_q-v\partial_q(u)\partial_p,
\eqno(1.15)$$
for $u,v\in A,\ p,q\in\ol{1,\ell},$
which is called a {\it Lie Algebra of generalized Witt type} [11,18,21].
\par
Let $\F[\DD]$ be the polynomial algebra with basis
$\{\ptl^{\mu}=\prod_{p=1}^\ell\ptl_p^{\mu_p}\,|\,\mu\in\Z_+^\ell\}$
(where $\ptl^\mu=1$ if $\mu=0$).
The $\F$-vector space
$$
\matrix{
\WW\!\!\!\!&=
\WW(\ell_1,\ell_2,\ell_3,\ell_4,\G)
\vs{4pt}\hfill\cr&
=\AA\otimes\F[\DD]=
{\rm span}
\{x^{\a,\vec i}\ptl^\mu\,|\,\a\in\G,\vec i\in\vec J,\mu\in\Z_+^\ell\},
\hfill\cr}
\eqno(1.16)$$
forms an associative algebra under the algebraic operation
$$
u\ptl^\mu\odot v\ptl^\nu=
u\sum_{\l\in \sZ_+^\ell}({}^{\dis\mu}_{\dis\l})\ptl^{\l}(v)\ptl^{\mu+\nu-\l},
\vsp
\eqno(1.17)$$
\vsp
where
$$({}^{\dis\mu}_{\dis\l})=
\prod_{p=1}^\ell ({}^{\dis\mu_p}_{\dis\l_p}),
\ \ \
\ptl^{\l}(v)=\prod_{p=1}^\ell\ptl_p^{\l_p}(v),
\vsp
\eqno(1.18)$$
(cf.~(1.12), (1.13)) for $u,v\in\AA,\mu,\nu,\l\in\Z_+^\ell$,
which induces a Lie algebra structure under the usual Lie bracket, called
a {\it Lie Algebra of Weyl type} [19,20]. In this paper, we shall call $\WW$
a {\it Lie Algebra of generalized differential operators} since it is a
generalization of the Lie algebra of the differential operators considered
in [4,9].
Note that using [19, Theorem 1.1], as in the proof of
[20, Theorem 1.1], one can derive that $\WW/\F$ is a simple Lie algebra,
which is non-graded and non-linear in general.
\par
It is well known that $W(0,0,0,1,\Z)$
is the centerless Virasoro algebra
and its 2-cohomology space is 1-dimensional such that the
Virasoro algebra is its universal central extension [3].
The 2-cohomology space for general
$W(\ell_1,\ell_2,\ell_3,\ell_4,\G)$ was
determined in [18].
The 2-cohomology space of
the Lie algebra $\WW(0,0,0,1,\Z)$ of the differential operators,
posed as an open problem by Kac [4], was determined by Li [9]. In [14], the
author determined all the 2-cohomology spaces of $\WW(0,0,0,n,\Z^n)$ for
$n\ge2$.
In [10], Li and Wilson determined the
2-cohomology spaces of the Lie algebra of the derivations
of the formal Laurent polynomial ring $\C((t))$,
the Lie algebra of
the differential operators of $\C((t))$ and
the Lie algebras of
the differential operators of $\C((t))\otimes\C^n$.
\par
In this paper, we shall determined all
the 2-cohomology spaces of the Lie algebras
$\WW(\ell_1,\ell_2,\ell_3,\ell_4,\G)$.
In [20], Zhao and the author determined the isomorphism classes
of the Lie algebras
$\WW(0,\ell_2,0,\ell_4,\G)$. However,
for general $\WW(\ell_1,\ell_2,\ell_3,\ell_4,\G)$,
the determination of the isomorphism classes
is still an unsolved problem.
We hope that the computation of the 2-cohomology spaces of
$\WW(\ell_1,\ell_2,\ell_3,\ell_4,\G)$ may help to determine the
isomorphism classes of $\WW(\ell_1,\ell_2,\ell_3,\ell_4,\G)$.
Our main result of the paper is stated in Theorem 3.1.
\par\
\par
\cl{\S2. \ 2-Cocycles on $\WW(\ell_1,\ell_2,\ell_3,\ell_4,\G)$
}
\par
In this section, we shall consider 2-cocycles
on the Lie algebra $\WW=\WW(\ell_1,\ell_2,\ell_3,\ell_4,\G)$,
i.e., we shall determine 2-cohomology group $H^2(\WW,\F)$.
\par
Take a basis of $\WW$:
$$
B=\{x^{\a,\vec i}\partial^\mu\,|\,
(\a,\vec i,\mu)\in\G\times\vec J\times\Z_+^\ell\}.
\eqno(2.1)$$
We fix a $\tau\in\G$ such that
$$
\tau_p\ne0,
\eqno(2.2)$$
for all $p\in\ol{\ell_1+1,\ell}.$
For any $\mu=(\mu_1,...,\mu_\ell)\in\Z^\ell$, we
define the {\it level} of $\mu$ to be
$|\mu|=\sum_{p=1}^\ell\mu_p$.
Define a total order on $\Z^\ell$ by: for $\mu,\nu\in\Z^\ell$, $\mu<\nu\Rla$
$$
|\mu|<|\nu|
\mbox{ or }
|\mu|=|\nu|
\mbox{ and for the first $p$ with $\mu_p\ne\nu_p$, we have $\mu_p<\nu_p$}.
\vms{-3pt}\eqno(2.3)$$
\par
{\bf Lemma 2.1}.
{\it Let $\psi$ be a 2-cocycle on $\WW$. Then there exists
a 2-cocycle $\phi$ equivalent to $\psi$ such
\vms{-1pt}
that
$$
\matrix{
\phi(t^{1_{[p]}}\partial_p,x^{\a,\vec i}\partial^\mu)
=0
\hfill&\mbox{if }p\in\ol{1,\ell'_3},
\vs{8pt}\hfill\cr
\phi(\partial_p,x^{\a,\vec i}\partial^\mu)
=0
\hfill&\mbox{if }p\in\ol{1,\ell},
\hfill\cr}
\eqno\matrix{(2.4)\vs{8pt}\cr(2.5)\cr}
\vms{-1pt}$$
for $(\a,\vec i,\mu)\in\G\times\vec J\times\Z_+^\ell$.}
\par
{\it Proof.} 
Define a $\F$-linear function $f:\WW\rar\F$ as follows.
\par
For $x^{\a,\vec i}\partial^\mu\in B$ with $\a\ne0$, let
$q$ be the minimal index such that $\a_q\ne0$. Define
$f(x^{\a,\vec i}\partial^\mu)$ inductively on $|i_q|$ by
$$
f(x^{\a,\vec i}\partial^\mu)
\!=\!\left\{\!\matrix{
\a_q^{-1}(\psi(\partial_q,x^{\a,\vec i}\partial^\mu)
-i_qf(x^{\a,\vec i-1_{[q]}}\partial^\mu))
\!\!\!\!\hfill&\mbox{if }
i_q\!\ge\!0,
\vs{4pt}\hfill\cr
-(1\!+\!\mu_q)^{-1}(\psi(t^{1_{[q]}}\partial_q,x^{\a,\vec i}\partial^\mu)
\!-\!\a_qf(x^{\a,\vec i+1_{[q]}}\partial^\mu))
\!\!\!\!\hfill&\mbox{if }i_q\!=\!-1,
\vs{4pt}\hfill\cr
(i_q+1)^{-1}(\psi(\partial_q,x^{\a,\vec i+1_{[q]}}\partial^\mu)
-\a_qf(x^{\vec\a,\vec i+1_{[q]}}\partial^\mu))
\!\!\!\!\hfill&\mbox{if }i_q\!\le\!-2.
\hfill\cr}\right.
\eqno(2.6)$$
Note that
$[\ptl_q,x^{\a,\vec i}\ptl^\mu]=\a_q x^{\a,\vec i}\ptl^\mu
+i_q x^{\a,\vec i}\ptl^\mu$ and
$[t^1\ptl_q,x^{\a,\vec i}\ptl^\mu]=\a_q x^{\a,\vec i+1_{[q]}}\ptl^\mu
+(i_q-\mu_q)x^{\a,\vec i}\ptl^\mu$.
If $i_q\ge1$, then $\vec i-1_{[q]}\in\vec J$, and so
$x^{\vec\a,\vec i-1_{[q]}}\partial_p\in B$; if $i_q\le-1$, then
$q\in\ol{\ell'_2+1,\ell'_3}$ (cf. (1.8)),
and $t^{1_{[q]}}\partial_q\in\WW$ and
$x^{\a,\vec i+1_{[q]}}\partial^\mu\in B$,
thus the right-hand side of (2.6) makes sense in all cases.
\par
For $t^{\vec i}\partial^\mu\in B$ with $\vec i\ne0$,
let $r$ be the minimal index such that $i_r\ne0$. Define
$$
f(t^{\vec i}\partial^\mu)
=\left\{\matrix{
(i_r-\mu_r)^{-1}\psi(t^{1_{[r]}}\partial_r,t^{\vec i}\partial^\mu)
\hfill&\mbox{if }i_r\ne \mu_r,
\vs{4pt}\hfill\cr
(i_r+1)^{-1}\psi(\partial_r,t^{\vec i+1_{[r]}}\partial^\mu)
\hfill&\mbox{if }i_r=\mu_r\mbox{ (and so $i_r\ge1$)}.
\hfill\cr}\right.
\eqno(2.7)$$
Note that in (2.7), since $i_r\ne0$, we have $r\le\ell'_3$, and so
all elements appearing in the right-hand side are in $\WW$.
\par
Finally, let $\ptl^\mu\in B$ with $\mu\in\Z_+^\ell$.
If $\ell'_3\ne0$, we define
$$
f(\partial^\mu)=
\psi(\partial_1,t^{1_{[1]}}\partial^\mu).
\vms{-4pt}\eqno(2.8)$$
On the other hand, if $\ell'_3=0$,
by using the
\vms{-10pt}formula
$$[x^\tau,x^{-\tau}\partial^{\mu+1_{[\ell]}}]
=-(\mu_\ell+1)\partial^\mu-
\sum_{0,1_{[\ell]\,}\ne\l\in\sZ_+^\ell}
(^{\dis \mu+1_{[\ell]}}_{\dis\,\ \ \l})
\prod_{p=1}^\ell\tau_p^{\l_p}\ptl^{\mu+1_{[\ell]}-\l},
\vms{-4pt}\eqno(2.9)$$
(cf.~(1.17), (1.18)),
we define $f(\partial^\mu)$ by induction on $\mu$ with respect to the order
defined in \vms{-10pt}(2.3):
$$
\matrix{\dis
f(\partial^\mu)=
-(\mu_\ell\!\!\!\!\!&+1)^{-1}(
\psi(x^\tau,x^{-\tau}\partial^{\mu+1_{[\ell]}})
\vs{4pt}\hfill\cr&\dis
+\sum_{0,1_{[\ell]\,}\ne\l\in\sZ_+^\ell}
(^{\dis \mu+1_{[\ell]}}_{\dis\,\ \ \l})
\prod_{p=1}^\ell\tau_p^{\l_p}f(\ptl^{\mu+1_{[\ell]}-\l})),
\hfill\cr}
\eqno(2.10)$$
where, by ordering (2.3),
we have $\mu+1_{[\ell]}-\l<\mu$ for 
$\l\in\Z_+^\ell\bs\{0,1_{[\ell]}\}$.
\par
Now set $\phi=\psi-\psi_f$ (cf.~(1.3)). Then by
(1.17), (2.6-8), (2.10), we have
$$
\matrix{
\phi(\partial_q,x^{\a,\vec i}\partial^\mu)=0\hfill&\mbox{if }\a\ne0,
\vs{8pt}\hfill\cr
\phi(t^{1_{[q]}}\partial_q,x^{\a,\vec i}\partial^\mu)=0
\hfill&\mbox{if }\a\ne0,\,i_q=-1,
\vs{8pt}\hfill\cr
\phi(t^{1_{[r]}}\partial_r,t^{\vec i}\partial^\mu)=0
\hfill&\mbox{if }i_r\ne\mu_r,
\vs{8pt}\hfill\cr
\phi(\partial_r,t^{\vec i}\partial^\mu)=0
\hfill&\mbox{if }i_r\ge2,\mu_r=i_r-1\mbox{ or }r=1,\vec i=1_{[1]},
\vs{8pt}\hfill\cr
\phi(x^{\tau},x^{-\tau}\partial^{\mu+1_{[\ell]}})=0
\hfill&
\mbox{if }\ell'_3=0,
\hfill\cr}
\eqno\matrix{(2.11)\vs{8pt}\cr
(2.12)\vs{8pt}\cr(2.13)\vs{8pt}\cr(2.14)\vs{8pt}\cr(2.15)\cr}
$$
where
$$
q={\rm min}\{q\in\ol{\ell_1+1,\ell}\,|\,\a_q\ne0\},\ \
r={\rm min}\{r\in\ol{1,\ell'_3}\,|\,i_r\ne0\}.
\eqno(2.16)$$
\par
We prove (2.4), (2.5) in three cases.
\par
{\it Case 1}: $\a\ne0$. 
\par
Let $q$ be as in (2.16). By (2.11), (1.3), (1.17), we obtain
$$
\matrix{
0\!\!\!\!&=\phi(\partial_q,\a_px^{\a,\vec i+1_{[p]}}\partial^\mu
+(i_p-\mu_p)x^{\a,\vec i}\partial^\mu)
\vs{4pt}\hfill\cr
&
=\phi(\partial_q,[t^{1_{[p]}}\partial_p,x^{\a,\vec i}\partial^\mu])
\vs{4pt}
\hfill\cr&
=\d_{q,p}\phi(\partial_q,x^{\a,\vec i}\partial^\mu)+
\phi(t^{1_{[p]}}\partial_p,\a_qx^{\a,\vec i}\partial^\mu
+i_qx^{\a,\vec i-1_{[q]}}\partial^\mu)
\vs{4pt}\hfill\cr
&
=\a_q\phi(t^{1_{[p]}}\partial_p,x^{\a,\vec i}\partial^\mu)
+i_q\phi(t^{1_{[p]}}\partial_p,x^{\a,\vec i-1_{[q]}}\partial^\mu),
\hfill\cr}
\eqno(2.17)$$
for $p\in\ol{1,\ell'_3},$ and
$$
\matrix{
0\!\!\!\!&=\phi(\partial_q,\a_px^{\a,\vec i}\partial^\mu
+i_px^{\a,\vec i-1_{[p]}}\partial^\mu)
\vs{4pt}\hfill\cr
&
=\phi(\partial_q,[\partial_p,x^{\a,\vec i}\partial^\mu])
\vs{4pt}\hfill\cr
&
=\a_q\phi(\partial_p,x^{\a,\vec i}\partial^\mu)
+i_q\phi(\partial_p,x^{\a,\vec i-1_{[q]}}\partial^\mu),
\hfill\cr}
\eqno(2.18)$$
for $p\in\ol{1,\ell}.$
If $i_q\ge0$, using (2.17), (2.18) and induction on $i_q$, we obtain (2.4),
(2.5). So assume that $i_q<0$.
First suppose $i_q=-1$. Then $q\in\ol{\ell'_2+1,\ell'_3}$ (cf.~(1.8)).
If $p=q$, then by (2.11), (2.12), we have (2.4), (2.5).
Assume that $p\ne q$.
Using (2.12), we have
$$
\matrix{
0\!\!\!\!&
=\phi(t^{1_{[q]}}\partial_q,[t^{1_{[p]}}\partial_p,
x^{\a,\vec i}\partial^\mu])
\vs{4pt}\hfill\cr
&=\a_q\phi(t^{1_{[p]}}\partial_p,x^{\a,\vec i+1_{[q]}}\ptl^\mu)
-(1+\mu_q)\phi(t^{1_{[p]}}\partial_p,x^{\a,\vec i}\partial^\mu)
\vs{4pt}\hfill\cr
&=-(1+\mu_q)\phi(t^{1_{[p]}}\partial_p,x^{\a,\vec i}\partial^\mu),
\hfill\cr}
\vms{-2pt}
\eqno(2.19)$$
for $p\in\ol{1,\ell'_3},$ where the last equality follows from that
$(\vec i+1_{[q]})_q=i_q+1=0$ and
\vms{-1pt}thus
$\phi(t^{1_{[p]}}\partial_p,x^{\a,\vec i+1_{[q]}}\ptl^\mu)=0$
by the argument above.
Thus (2.4) holds if $i_q=-1$. Then
using (2.17) and induction on $|i_q|$,
we see that (2.4) holds for all $i_q$.
For (2.5), if $i_q=-1$, we have
$$
\matrix{
0\!\!\!\!&
=\phi(t^{1_{[q]}}\partial_q,[\partial_p,x^{\a,\vec i}\partial^\mu])
\vs{4pt}\hfill\cr
&=-\d_{p,q}\phi(\partial_p,x^{\a,\vec i}\partial^\mu)
+\a_q\phi(\partial_p,x^{\a,\vec i+1_{[q]}}\partial^\mu)
-(1+\mu_q)\phi(\partial_p,x^{\a,\vec i}\partial^\mu)
\vs{4pt}\hfill\cr
&=-(1+\d_{p,q}+\mu_q)\phi(\partial_p,x^{\a,\vec i}\partial^\mu),
\hfill\cr}
\vms{-2pt}
\eqno(2.20)$$
for $p\in\ol{1,\ell}$, where the last equality follows again from
that $(\vec i+1_{[q]})_q=0$.
Thus (2.5) holds if $i_q=-1$.
Then by (2.18) and induction on $|i_q|$, we deduce (2.5) for all $i_q$.
\par
{\it Case 2}: $\a=0,\vec i=0$.
\vms{-1pt}\par
First, consider (2.5).
If $p\in\ol{1,\ell'_3}$, then
$$
0=\phi(t^{1_{[p]}}\ptl_p,[\ptl_p,\ptl^\mu])
=-(1+\mu_p)\phi(\ptl_p,\ptl^\mu),
\eqno(2.21)$$
and we have (2.5) in this case. Thus assume that $p\in\ol{\ell'_3+1,\ell}$.
If $\ell'_3\ge1$, then
$$
\phi(\ptl_p,\ptl^\mu)
=\phi(\ptl_p,[\ptl_1,t^{1_{[1]}}\ptl^\mu])
=\phi(\ptl_1,[\ptl_p,t^{1_{[1]}}\ptl^\mu])=0.
\eqno(2.22)$$
On the other hand, if $\ell'_3=0$, then by
\vsp
(2.15),
$$
\matrix{
0\!\!\!\!&=\tau_p\phi(x^{\tau},x^{-\tau}\ptl^{\mu+1_{[\ell]}})
\vs{4pt}\hfill\cr&\dis
=\phi([\ptl_p,x^{\tau}],x^{-\tau}\ptl^{\mu+1_{[\ell]}})
=-\sum_{0\ne\l\in\sZ_+^\ell}
(^{\dis \mu+1_{[\ell]}}_{\dis\,\ \ \l})
\prod_{q=1}^\ell\!\tau_q^{\l_q}
\phi(\ptl_p,\ptl^{\mu+1_{[\ell]}-\l}),
\hfill\cr}
\vsp
\eqno(2.23)$$
where the last equality follows from that
$\phi([\ptl_p,x^{-\tau}\ptl^{\mu+1_{[\ell]}}],x^\tau)=0$ by (2.15).
Using induction on $\mu$ in (2.23) gives $\phi(\ptl_p,\ptl^\mu)=0$.
Next, consider (2.4). So suppose $p\in\ol{1,\ell'_3}$. Then $\ell'_3\ge1$.
If $p=1$, we have
$$
\phi(t^{1_{[1]}}\ptl_1,\ptl^\mu)
=2^{-1}\phi([\ptl_1,t^{2_{[1]}}\ptl_1],\ptl^\mu)
=2^{-1}\phi(\ptl_1,[t^{2_{[1]}}\ptl_1,\ptl^\mu])
=0,
\eqno(2.24)$$
where the last equality follows
from (2.14), (2.22) and that
$$[t^{2_{[1]}}\ptl_1,\ptl^\mu]=-2\mu_1t^{1_{[1]}}\ptl^\mu-\mu_1(\mu_1-1)\ptl^{\mu-1_{[1]}}.
\eqno(2.25)$$
Thus assume that $p\in\ol{2,\ell'_3}$. If $\mu_1\ne0$, then by (2.24),
$$
0=\phi(t^{1_{[1]}}\ptl_1,[t^{1_{[p]}}\ptl_p,\ptl^\mu])
=\phi(t^{1_{[p]}}\ptl_p,[t^{1_{[1]}}\ptl_1,\ptl^\mu])
=-\mu_1\phi(t^{1_{[p]}}\ptl_p,\ptl^\mu).
\eqno(2.26)$$
On the other hand, if $\mu_1=0$, then
$$
\matrix{
\phi(t^{1_{[p]}}\ptl_p,\ptl^\mu)
\!\!\!\!&=\phi([\ptl_1,t^{1_{[1]}+1_{[p]}}\ptl_p],\ptl^\mu)
\vs{4pt}\hfill\cr&
=\phi(\ptl_1,[t^{1_{[1]}+1_{[p]}}\ptl_p,\ptl^\mu])
=-\mu_p\phi(\ptl_1,t^{1_{[1]}}\ptl^\mu)=0,
\hfill\cr}
\eqno(2.27)$$
where the last equality follows from (2.14).
\par
{\it Case 3}: $\a=0,\vec i\ne0$.
\vs{-0.1pt}\par
Since $\vec i\ne0$, we have $\ell'_3\ge1$ (cf.~(1.8)).
Let $r$ be as in (2.16). First consider (2.4).
If $i_r\ne\mu_r$, then by (2.13), we have
$$
\matrix{
0\!\!\!\!&=\phi(t^{1_{[r]}}\partial_r,
[t^{1_{[p]}}\partial_p,t^{\vec i}\partial^\mu])
\vs{4pt}\hfill\cr&
=\phi(t^{1_{[p]}}\partial_p,
[t^{1_{[r]}}\partial_r,t^{\vec i}\partial^\mu])
=(i_r-\mu_r)\phi(t^{1_{[p]}}\partial_p,t^{\vec i}\partial^\mu),
\hfill\cr}
\eqno(2.28)$$
thus we have (2.4) in this case.
So assume that $i_r=\mu_r$ (and so $i_r\ge0$). Then
$$
\matrix{
\phi(t^{1_{[p]}}\partial_p,t^{\vec i}\partial^\mu)
\!\!\!\!&
=(i_r+1)^{-1}\phi(t^{1_{[p]}}\partial_p,
[\partial_r,t^{\vec i+1_{[r]}}\partial^\mu])
\vs{4pt}\hfill\cr&
=(i_r+1)^{-1}(-\d_{r,p}+(i_p+\d_{p,r}-\mu_p))\phi(\partial_r,
t^{\vec i+1_{[r]}}\partial^\mu)
\vs{4pt}\hfill\cr
&=0,
\hfill\cr}
\eqno(2.29)$$
where the last equality follows from (2.13).
This proves (2.4) in this case.
\par
Next consider (2.5).
We can write
$$
t^{\vec i}\partial^\mu=
\left\{\matrix{
(i_1+1)^{-1}[\partial_1,t^{\vec i+1_{[1]}}\partial^\mu]
\hfill&\mbox{if }i_1\ne-1,
\vs{4pt}\hfill\cr
-(1+\mu_1)^{-1}[t^{1_{[1]}}\partial_1,t^{\vec i}\partial^\mu]
\hfill&\mbox{if }i_1=-1.
\hfill\cr}\right.
\eqno(2.30)$$
Thus by (1.2) and (2.4), we have
$$
\phi(\partial_p,t^{\vec i}\partial^\mu)=
\left\{\matrix{
(i_1+1)^{-1}(i_p+\d_{p,1})\phi(\partial_1,t^{\vec i+1_{[1]}-1_{[p]}}
\partial^\mu)
\hfill&\mbox{if }i_1\ne-1,
\vs{4pt}\hfill\cr
-(1+\mu_1)^{-1}\d_{p,1}\phi(\partial_1,t^{\vec i}\partial^\mu)
\hfill&\mbox{if }i_1=-1.
\hfill\cr}\right.
\eqno(2.31)$$
Thus it suffices to prove (2.5) for $p=1$.
Using (2.4), we have
$$
0
=\phi(t^{1_{[1]}}\partial_1,[\partial_1,t^{\vec i}\partial^\mu])
=(-1+i_1-\mu_1)\phi(\partial_1,t^{\vec i}\partial^\mu).
\eqno(2.32)$$
By (2.32), it remains to consider the case $i_1=\mu_1+1$.
If $\mu_1\ge1$ or $\vec i=1_{[1]}$, the result follows from (2.14).
Thus assume $\mu_1=0$ and $\vec i\ne1_{[1]}$. So $i_{r'}\ne0$ for some
$r'\in\ol{2,\ell'_3}$. Let $r'\in\ol{2,\ell'_3}$ be
the minimal with $i_{r'}\ne0$. If $i_{r'}\ne\mu_{r'}$, then by (2.4),
$$
0=\phi(t^{1_{[r']}}\ptl_{r'},[\ptl_1,t^{\vec i}\ptl^\mu])
=(i_{r'}-\mu_{r'})\phi(\ptl_1,t^{\vec i}\ptl^\mu).
\eqno(2.33)$$
On the other hand, if $i_{r'}=\mu_{r'}$ (and so $i_{r'}\ge0$), then
$$
\matrix{
\phi(\ptl_1,t^{\vec i}\ptl^\mu)
\!\!\!\!&=(i_{r'}+1)^{-1}\phi(\ptl_1,[\ptl_{r'},t^{\vec i+1_{[r']}}\ptl^\mu])
\vs{4pt}\hfill\cr&
=(i_{r'}+1)^{-1}\phi(\ptl_{r'},t^{\vec i-1_{[1]}+1_{[r']}}\ptl^\mu)
=0,
\hfill\cr}
\eqno(2.34)$$
where the last equality follows from (2.14). This proves Lemma 2.1.
\qed\par
$\rm L\sc EMMA$ 2.2.
{\it If $\ell'_2=\ell_1+\ell_2\ge1$, then $H^2(\WW,\F)=0$.}
\par
{\it Proof.} By Lemma 2.1, suppose $\phi$ is 2-cocycle on $\WW$ satisfying
(2.4), (2.5). We shall prove
$$
\phi(x^{\a,\vec i}\partial^\mu,x^{\b,\vec j}\partial^\nu)=0,
\eqno(2.35)$$
for all $x^{\a,\vec i}\partial^\mu,x^{\b,\vec j}\partial^\nu\in B.$
By (2.5), we obtain
$$
\matrix{
0=\!\!\!\!&
\phi(\partial_1,
[x^{\a,\vec i}\partial\mu,x^{\b,\vec j}\partial^\nu])
\vs{4pt}\hfill\cr\hfill=\!\!\!\!&
(\a_1+\b_1)\phi(x^{\a,\vec i}\partial^\mu,
x^{\b,\vec j}\partial^\mu)
\vs{4pt}\hfill\cr&
\ \ \ \,\
+\,i_1\phi(x^{\a,\vec i-1_{[1]}}\partial^\mu,x^{\b,\vec j}\partial^\nu)
+j_1\phi(x^{\a,\vec i}\partial^\mu,x^{\b,\vec j-1_{[1]}}\partial^\nu).
\hfill\cr}
\eqno(2.36)$$
Since $\ell'_2\ge1$, we have $i_1,j_1\ge0$ (cf.~(1.8)). If
$\a_1+\b_1\ne0$, then using induction on $i_1+j_1$, we immediately get
(2.35); otherwise, by substituting $\vec i$ by $\vec i+1_{[1]}$ in (2.36),
and using induction on $j_1\ge0$, we again deduce (2.35).
\qed
\par\ \par
\cl{\bf 3. The main results}
\par
Before given the main result of this paper, for convenience, we introduce the
following notations
$$
[\a]_\mu=\a(\a-1)\cdots(\a-\mu+1),\ \
(^{\dis\a}_{\dis\mu})=(\mu!)^{-1}[\a]_\mu,
\eqno(3.1)$$
for $\a\in\F,\mu\in\Z_+$. We shall also use the following conventions
$$
\matrix{
\a|_{\a=0}=\rb{-5pt}{\mbox{$^{\dis\rm lim}_{\a\rar0}$}}\a,
\vs{4pt}\hfill\cr
\dis{1\over k!}=0\mbox{ if }k\le-1,
\hfill\cr}
\eqno\matrix{(3.2)\vs{6pt}\hfill\cr(3.3)\cr}$$
where (3.2) shall be interpreted as follows: whenever we take a value
$\a$ to be zero in an expression, we shall do it by taking the limit
$\a\rar0$.
\par
{\bf Theorem 3.1}.
{\it (1) If $\ell_1+\ell_2\ne0$ or $\ell\ge2$, then $H^2(\WW,\F)=0$.
\par
(2) If $\ell=\ell_4=1$, then $H^2(\WW,\F)=\F\ol\phi_0$, where $\ol\phi_0$ is the
cohomology class of $\phi_0$ defined by
$$\phi_0(x^\a\ptl^\mu,x^\b\ptl^\nu)
=\d_{\a+\b,0}(-1)^\mu\mu!\nu!(^{\ \dis\ \a+\mu}_{\dis \mu+\nu+1}),
\eqno(3.4)$$
for all $\a,\b\in\G\subseteq\F,\,\mu,\nu\in\Z_+,$ where $\ptl=\ptl_1$.
\par
(3) Suppose $\ell=\ell_3=1$.
Then for any $\g\in\G$, there exists a cohomology class
$\ol\phi_\g$ defined
\vsp
by
$$
\matrix{\dis
\phi_\g(x^{\a,i}\ptl^\mu,x^{\b,j}\ptl^\nu)
\vs{4pt}\hfill\cr\dis\ \ \ \ \
=\d_{\a+\b,\g}(-1)^\mu\mu!\nu!
\sum_{s=0}^{\mu+\nu+1}
(^{\ssc\,\dis i}_{\dis s}){\a^{\mu+\nu+1-s}\over(\mu+\nu+1-s)!}
\cdot{\g^{s-i-j-1}\over(s-i-j-1)!}\ ,
\hfill\cr}
\vsp
\eqno(3.5)$$
for all $(\a,i,\mu),(\b,j,\nu)\in\G\times\Z\times\Z_+$
(if $\g=0$, then by conventions (3.2), (3.3),
${\g^{s-i-j-1}\over(s-i-j-1)!}=\d_{i+j,s-1}$).
Furthermore,
$$
H^2(\WW,\F)=\prod_{\a\in\G}\F\ol\phi_\g,
\eqno(3.6)$$
is the direct product of all $\F\ol\phi_\g,\g\in\G.$
}
\par
{\it Proof.} (1) follows from Lemma 2.2 if $\ell'_2\ge1$.
Next we suppose that $\ell'_2=0$
and that $\phi$ is a 2-cocycle on $\WW$ satisfying (2.4-5),
(2.11-15).
We shall divide the discussion into three cases.
\par
{\it Case 1}: $\ell=\ell_4=1$.
\par
Then $\WW=\WW(0,0,0,1,\G)={\rm span}\{x^\a\ptl^\mu\,|\,\a\in\G,\mu\in\Z_+\}$
has the commutation relations (cf.~(1.17))
$$
[x^\a\ptl^\mu,x^\b\ptl^\nu]=\sum_{s=0}^{\mu+\nu}
((^{\dis\mu}_{\ssc\,\dis s})\b^s-
(^{\dis\nu}_{\ssc\,\dis s})\a^s)x^{\a+\b}\ptl^{\mu+\nu-s}.
\eqno(3.7)$$
Note that when $\G=\Z$, (3.4) is precisely the 2-cocycle
defined in [4] and [9], thus as in [4], one can verify that,
in general, (3.4) defines a nontrivial 2-cocycle on $\WW(0,0,0,1,\G)$.
\par
If necessary by replacing $\G$ by $\tau^{-1}\G$, we can suppose that the
element $\tau$ of $\G$ defined in (2.2) is $\tau=1$.
Now let $\phi$ be a 2-cocycle satisfying (2.5), (2.15).
Let $\phi_1=\phi-\phi(x^1,x^{-1})\phi_0$.
Then $\phi_1$ is a 2-cocycle satisfying (2.5), (2.15) and
$$
\phi_1(x^1,x^{-1})=0.
\eqno(3.8)$$
We shall prove $\phi_1=0$. We have
$$
\matrix{
(\a+\b)\phi_1(x^\a\ptl^\mu,x^\b\ptl^\nu)=
\phi_1(\ptl,[x^\a\ptl^\mu,x^\b\ptl^\nu])=0,
\vs{8pt}\hfill\cr
\phi_1(x^\a,x^{-\a})=
\phi_1([x^{\a-1}\ptl,x^1],x^{-\a})=
\a\phi_1(x^1,x^{-1})=0,
\hfill\cr}
\eqno\matrix{(3.9)\vs{8pt}\cr(3.10)\cr}$$
for $\a,\b\in\G,\mu,\nu\in\Z_+$. Thus
$\phi_1(x^\a\ptl^\mu,x^\b\ptl^\nu)=0$ for $\a+\b\ne0$ or $\mu+\nu=0$.
Using induction on $\mu+\nu$, suppose that for $n\ge0$ we have proved
$$
\phi_1(x^\a\ptl^\mu,x^\b\ptl^\nu)=0,
\eqno(3.11)$$
for all $\a,\b\in\G,\,\mu,\nu\in\Z_+$ with $\mu+\nu\le n$.
Then by (3.7), (3.11),
$$
\matrix{
\phi_1(x^\a\ptl^{n+1},x^{-\a})
\!\!\!\!&=
(n+2)^{-1}\phi_1([x^{\a-1}\ptl^{n+2},x^1],x^{-\a})
\vs{4pt}\hfill\cr&
=\a\phi_1(x^1,x^{-1}\ptl^{n+1})
\vs{4pt}\hfill\cr&
=0,
\hfill\cr}
\eqno(3.12)$$
where the last equality follows from (2.15); and if $1\le\mu\le n$,
$$
\matrix{
\phi_1(x^\a\ptl^\mu,x^{-\a}\ptl^{n+1-\mu})
\!\!\!\!&=-((n+2-\mu)\a)^{-1}\phi_1(x^\a\ptl^\mu,[\ptl^{n+2-\mu},x^{-\a}])
\vs{4pt}\hfill\cr&
=(n+2-\mu)^{-1}\mu\phi_1(\ptl^{n+2-\mu},\ptl^{\mu-1})
\vs{4pt}\hfill\cr&
=0,
\hfill\cr}
\eqno(3.13)$$
where the second equality follows from (3.12) and
the last equality is obtained by induction on $\mu$.
This proves that (3.11) holds for $\mu+\nu=n+1$. Hence $\phi_1=0$.
\par
{\it Case 2}: $\ell=\ell_3=1$.
\par
Then $\vec J=\Z$ and $\vec i=i$, and
$\WW=\WW(0,0,1,0,\G)={\rm span}\{x^{\a,i}\ptl^\mu\,|\,
(\a,i,\mu)\in\G\times\Z\times\Z_+\}$
has the commutation relations (cf.~(1.17),
\vsp
(3.7))
$$
\matrix{
[x^{\a,i}\ptl^\mu,\!\!\!\!\!&\dis x^{\b,j}\ptl^\nu]
=\sum_{s=0}^{\mu+\nu}
((^{\dis\mu}_{\ssc\,\dis s})x^{\a,i}\ptl^s(x^{\b,j})
-(^{\dis\nu}_{\ssc\,\dis s})\ptl^s(x^{\a,i})x^{\b,j})\ptl^{\mu+\nu-s}
\vs{4pt}\hfill\cr&
\dis=\sum_{s=0}^{\mu+\nu}\sum_{r=0}^s
((^{\dis\mu}_{\ssc\,\dis s})(^{\dis s}_{\dis r})
[j]_r\b^{s-r}
-(^{\dis\nu}_{\ssc\,\dis s})(^{\dis s}_{\dis r})
[i]_r\a^{s-r})x^{\a+\b,i+j-r}\ptl^{\mu+\nu-s}.
\hfill\cr}
\vsp
\eqno(3.14)$$
\par
Let $\phi$ be a 2-cocycle satisfying (2.4), (2.5). Note that by
(1.2), one can derive that $\phi(1,v)=0$ for all $v\in\WW$.
Then by (2.4), (2.5), we have
$$
\matrix{
0\!\!\!\!&=\phi(\ptl,[x^{\a,i}\ptl^\mu,t^1])
=\phi([\ptl,x^{\a,i}\ptl^\mu],t^1])+\phi(x^{\a,i}\ptl^\mu,[\ptl,t^1])
\vs{4pt}\hfill\cr&
=\phi(\ptl(x^{\a,i})\ptl^\mu,t^1])
=\a\phi(x^{\a,i}\ptl^\mu,t^1)+i\phi(x^{\a,i-1}\ptl^\mu,t^1),
\vs{8pt}\hfill\cr
0\!\!\!\!&=\phi(t^1\ptl,[x^{\a,i-1}\ptl^\mu,t^1])
=\a\phi(x^{\a,i}\ptl^\mu,t^1)+(i-\mu)\phi(x^{\a,i-1}\ptl^\mu,t^1).
\hfill\cr}
\eqno\matrix{\cr(3.15)\vs{8pt}\cr(3.16)\cr}$$
Comparing the coefficients of both terms of the right-hand sides of
the above two equations, we deduce
$$
\matrix{
\phi(x^{\a,i}\ptl^\mu,t^1)=0
\mbox{ if }\mu\ne0,\mbox{ \ or }\a\ne0,i\ge0,\mbox{ \ or }\a=0,i\ne-1,
\vs{4pt}\hfill\cr
\dis\phi(x^{\a,i},t^1)=-{\a^{-i-1}\over(-i-1)!}b_\a\mbox{ if }
\a\ne0\mbox{ and }i\le-1,
\hfill\cr}
\eqno\matrix{(3.17)\vs{8pt}\cr(3.18)\cr}
\vms{-4pt}
$$
where $b_\a=\phi(t^1,x^{\a,-1})\in\F$ for $\a\in\G$.
Using convention (3.2), (3.3), we
\vsp
have
$$
\phi(x^{\a,i},t^1)=-{\a^{-i-1}\over(-i-1)!}b_\a\mbox{ for all }
(\a,i)\in\G\times\Z.
\vsp
\eqno(3.19)$$
\vsp
Then
$$
\matrix{
\phi(x^{\a,i},x^{\b,j}\ptl^\nu)
\!\!\!\!&
\dis={1\over\nu+1}
\phi(x^{\a,i},[x^{\b,j}\ptl^{\nu+1},t^1])
\vs{4pt}\hfill\cr&
\dis=-{1\over\nu+1}
\phi(\ptl^{\nu+1}(x^{\a,i})x^{\b,j},t^1])
\vs{4pt}\hfill\cr&
\dis={1\over\nu+1}
\sum_{s=0}^{\nu+1}(^{\dis\nu+1}_{\dis\ \ \ s})[i]_s
\a^{\nu+1-s}\phi(t^1,x^{\a+\b,i+j-s})
\vs{4pt}\hfill\cr&
\dis=
\nu!\sum_{s=0}^{\nu+1}
(^{\ssc\,\dis i}_{\dis s}){\a^{\nu+1-s}\over(\nu+1-s)!}
\cdot{\g^{s-i-j-1}\over(s-i-j-1)!}b_\g,
\hfill\cr}
\vsp
\eqno(3.20)$$
where $\g=\a+\b$, and where, the second equality follows from (3.17)
and the last equality follows from (3.19);
\vsp
and if $\mu>0$,
$$
\matrix{
\phi(x^{\a,i}\ptl^\mu,x^{\b,j}\ptl^\nu)
\!\!\!\!&
\dis={1\over\nu+1}
\phi(x^{\a,i}\ptl^\mu,[x^{\b,j}\ptl^{\nu+1},t^1])
\vs{4pt}\hfill\cr&
\dis={\mu\over\nu+1}
\phi(x^{\b,j}\ptl^{\nu+1},x^{\a,i}\ptl^{\mu-1})
\vs{4pt}\hfill\cr&
\dis=(-1)^\mu{\mu!\nu!\over(\mu+\nu)!}
\phi(x^{\a,i},x^{\b,j}\ptl^{\mu+\nu})
\vs{6pt}\hfill\cr&
=b_r\times \mbox{the right-hand side of (3.5)},
\hfill\cr}
\vsp
\eqno(3.21)$$
where the second equality follows from (3.17) and the third equality
is obtained by induction on $\mu$ and the last equality follows from (3.20).
\par
First assume that $\phi=\psi_f$ is a 2-coboundary
defined by a linear function $f$ (cf.~(1.3)).
Then by (2.4), (2.5), we have
$$
\matrix{
\a f(x^{\a,i}\ptl^\mu)+(i-1-\mu)f(x^{\a,i-1}\ptl^\mu)
=\phi(t^1\ptl,x^{\a,i-1}\ptl^\mu)=0,
\vs{8pt}\hfill\cr
\a f(x^{\a,i}\ptl^\mu)+i f(x^{\a,i-1}\ptl^\mu)
=\phi(\ptl,x^{\a,i}\ptl^\mu)=0.
\hfill\cr}
\eqno\matrix{(3.22)\vs{8pt}\cr(3.23)\cr}$$
Thus $f=0$, and so $\phi=0$ and $b_\g=\phi(t^1,x^{\g,-1})=0$
for all $\g\in\G$.
Thus (3.20), (3.21) show
\vsp
that
$$
\phi=\prod_{\g\in\G}b_\g\phi_\g,
\vsp
\eqno(3.24)$$
is the direct product of $b_\g\phi_\g$ for $\g\in\G$, where
$\phi_\g$ is defined in (3.5), i.e., we have (3.6)
as long as each $\phi_\g$ is a 2-cocycle. It remains
to prove that $\phi$ defined in (3.24) is a 2-cocycle for any $\b_\g\in\F,
\g\in\G$.
\par
So now let $\varphi=\prod_{\g\in\G}b_\g\phi_\g$. It is straightforward
to check that formulas (3.15-21) hold for $\varphi$.
By (3.21), the second equality of (3.20), we
\vsp
have
$$
\matrix{
\varphi(x^{\a,i}\ptl^\mu,x^{\b,j}\ptl^\nu)
\!\!\!\!&
\dis=(-1)^\mu{\mu!\nu!\over(\mu+\nu+1)!}
\varphi(t^1,\ptl^{\mu+\nu+1}(x^{\a,i})x^{\b,j})
\vs{4pt}\hfill\cr&
\dis=(-1)^{\nu+1}{\mu!\nu!\over(\mu+\nu+1)!}
\varphi(t^1,x^{\a,i}\ptl^{\mu+\nu+1}(x^{\b,j}))
\vs{8pt}\hfill\cr&
=-\varphi(x^{\b,j}\ptl^\nu,x^{\a,i}\ptl^\mu),
\hfill\cr}
\eqno(3.25)$$
where the second equality follows from the third equality of (3.15)
so that
$$
\varphi(t^1,\ptl(x^{\a,i})x^{\b,j})=
-\varphi(t^1,x^{\a,i}\ptl(x^{\b,j})),
\eqno(3.26)$$
and the last equality of (3.25) follows from the first equality. This
proves the
skew-symmetry (1.1). To prove (1.2), recall that $\WW$ is an associative
algebra under algebraic operation (1.17).
By (1.17) (cf.~(3.14)), and the second equality of (3.25),
\vsp
we have
$$
\matrix{\dis
\varphi(u\ptl^\mu,\!\!\!\!\!&v\ptl^\nu\odot w\ptl^\l)
\vs{4pt}\hfill\cr&\dis=
\sum_{s=0}^\nu(-1)^{\nu+\l+1-s}
{\mu!(\nu+\l-s)!\over(k+1-s)!}
(^{\dis\nu}_{\ssc\,\dis s})
\varphi(t^1,u\ptl^{k+1-s}(v\ptl^s(w))),
\hfill\cr}
\vsp
\eqno(3.27)$$
for $u,v,w\in\AA$, where $k=\mu+\nu+\l$. Using shifted version of (3.27) and
(3.26),
\vsp
we have
$$
\matrix{\dis
\varphi(v\ptl^\nu,\!\!\!\!\!&w\ptl^\l\odot u\ptl^\mu)
\vs{4pt}\hfill\cr&\dis=
\sum_{s=0}^\l(-1)^{\nu+s}
{\nu!(\l+\mu-s)!\over(k+1-s)!}(^{\dis\l}_{\ssc\,\dis s})
\varphi(t^1,u\ptl^s(\ptl^{k+1-s}(v)w)),
\vs{8pt}\hfill\cr\dis
\varphi(w\ptl^\l,\!\!\!\!\!&u\ptl^\mu\odot v\ptl^\nu)
\vs{4pt}\hfill\cr&\dis=
\sum_{s=0}^\mu(-1)^\l
{\l!(\nu+\mu-s)!\over(k+1-s)!}(^{\dis\mu}_{\ssc\,\dis s})
\varphi(t^1,u\ptl^s(v)\ptl^{k+1-s}(w)).
\hfill\cr}
\eqno\matrix{\cr(3.28)\cr\cr\vs{15pt}\cr(3.29)\cr}
\vsp
$$
Denote the right-hand sides of (3.27), (3.28) and (3.29) by
$$
\sum_{s=0}^{k+1}a_{p,s}\phi(t^1,u\ptl^s(v)\ptl^{k+1-s}(w))
\mbox{ for }p=1,2\mbox{ and 3 respectively.}
$$
Using
$(1+x)^{k+1-s}(1+x)^{-(\l+1)}=(1+x)^{\mu+\nu-s}$, we deduce the following
binomial
\vsp
formula
$$
\sum_q(-1)^q(^{\dis k+1-s}_{\dis\ \ \nu-q})(^{\dis\l+q}_{\,\dis\ \ q})=
(^{\dis\mu+\nu-s}_{\,\dis\ \ \ \ \,\nu}).
\vsp
\eqno(3.30)$$
Using (3.30), we can deduce that if $s\le\mu$,
\vsp
then
$$
\matrix{
a_{1,s}\!\!\!\!&\dis=\sum_{q=0}^\nu
(-1)^{\nu+\l+1-q}
{\mu!(\nu\!+\!\l\!-\!q)!\over(k+1-q)!}
(^{\dis\nu}_{\ssc\,\dis q})
(^{\dis k\!+\!1\!-\!q}_{\dis\!\!\ \ \ \ \ s})
\vs{4pt}\hfill\cr&\dis
=(-1)^{\l+1}{\l!(\nu\!+\!\mu\!-\!s)!\over(k+1-s)!}(^{\dis\mu}_{\ssc\,\dis s})
\vs{6pt}\hfill\cr&\dis
=-a_{3,s},
\hfill\cr}
\vsp
\eqno(3.31)$$
and $a_{2,s}=0$; and if $\mu<s\le\mu+\nu$, then $a_{1,s}=a_{2,s}=a_{3,s}=0$;
and if $\mu+\nu<s\le k+1$,
\vsp
then
$$
\matrix{
a_{1,s}\!\!\!\!&\dis
=\sum_{q=0}^\nu
(-1)^{\nu+\l+1-q}
{\mu!(\nu+\l-q)!\over(k+1-q)!}
(^{\dis\nu}_{\ssc\,\dis q})
(^{\dis k+1-q}_{\dis\ \ \ \ \ s})
\vs{4pt}\hfill\cr&\dis
=
\sum_{q=0}^\l(-1)^{\nu+q+1}
{\nu!(\l+\mu-q)!\over(k+1-q)!}
(^{\dis\l}_{\ssc\,\dis q})
(^{\dis\ \ \ \ \ q}_{\dis k+1-s})
\vs{6pt}\hfill\cr&\dis
=-a_{2,s},
\hfill\cr}
\vsp
\eqno(3.32)$$
and $a_{3,s}=0$.
This proves that the sum of (3.27), (3.28) and (3.29) is zero.
In particular,
$\varphi$ satisfies (1.2). So, $\varphi$ is a 2-cocycle on $\WW$.
\par
{\it Case 3}: $\ell\ge2$.
\par
First suppose $\ell'_3\ge1$. As in (3.15-17), we have
$$
\phi(t^{1_{[1]}},x^{\a,\vec i}\ptl^\mu)=0
\eqno(3.33)$$
if $\mu_1\ne0$. Assume that $\mu_1=0$. We can
\vsp
write
$$
x^{\a,\vec i}\ptl^\mu
\equiv{1\over\tau_\ell(\mu_\ell+1)}
[x^{\a-\tau,\vec i}\ptl^{\mu+1_{[\ell]}},x^\tau]
\ ({\rm mod\,}\WW_{(\mu)}),
\vsp
\eqno(3.34)$$
where $\WW_{(\mu)}={\rm span}\{x^{\a,\vec i}\ptl^\nu\,|\,
(\a,\vec i,\nu)\in\G\times\vec J\times\Z_+^\ell,\nu<\mu\}$
(cf.~(2.3)). Then using induction on $\mu$, by (3.34), (1.2),
we can prove that
(3.33) also holds for $\mu_1=0$. Also as in (3.20), (3.21), we have
$$
\phi(x^{\a,\vec i}\ptl^\mu,x^{\b,\vec j}\ptl^\nu)
=(\nu_1+1)^{-1}
\phi([x^{\a,\vec i}\ptl^\mu,x^{\b,\vec j}\ptl^{\nu+1_{[1]}}],t^{1_{[1]}})=0
\mbox{ if }\mu_1=0,
\eqno(3.35)$$
\vsp
where the last equality follows from (3.33),
and
$$
\phi(x^{\a,\vec i}\ptl^\mu,x^{\b,\vec j}\ptl^\nu)
=-{\mu_1\over\nu_1+1}\phi(x^{\a,\vec i}\ptl^{\mu-1_{[1]}},x^{\b,\vec j}
\ptl^{\nu+1_{[1]}})=0
\mbox{ if }\mu_1\ne0.
\vsp
\eqno(3.36)$$
where the last equality is obtained by induction on $\mu_1$. Thus $\phi=0$.
\par
Next suppose $\ell'_3=0$. By replacing $\ptl_p$
\vsp
by
$$
\ptl'_p=\ptl_p-{\tau_p\over\tau_\ell}\ptl_\ell
\mbox{ for $p=1,...,\ell-1$,}
\vsp
\eqno(3.37)$$
we can suppose
$$\ptl_p(x^\tau)=0
\mbox{ for $p=1,...,\ell-1$.}
\eqno(3.38)$$
Then using (2.5), (2.15) and (1.2), we obtain
$$
\matrix{
\phi(x^\a\ptl^\mu,x^\b\ptl^\nu)=0\hfill&\mbox{if }\a+\b\ne0,
\vs{8pt}\hfill\cr
\phi(x^\tau,x^{-\tau}\ptl^\mu)=0\hfill&\mbox{if }\mu_\ell\ne0.
\hfill\cr
}
\eqno\matrix{(3.39)\vs{8pt}\cr(3.40)\cr}$$
Assume that $\mu_\ell=0$. Suppose we have proved that
$\phi(x^\tau,x^{-\tau}\ptl^\nu)=0$ for $\nu<\mu$. Choose
$\b\in\G$ with $\b_{\ell-1}\ne0$, then we
have
$$
\phi(x^\tau,x^{-\tau}\ptl^\mu)
=(\b_{\ell-1}(\mu_{\ell-1}+1))^{-1}
\phi(x^\tau,[x^{-\tau+\b}\ptl^{\mu+1_{[\ell-1]}},x^\b])
=0,
\eqno(3.41)$$
where the first equality follows from the assumption that
$\phi(x^\tau,x^{-\tau}\ptl^\nu)=0$ for $\nu<\mu$, and
the last equality follows from (1.2) and (3.38).
This proves that $\phi(x^\tau,x^{-\tau}\ptl^\mu)=0$ for $\mu\in\Z_+$.
For any $\a\in\G,\nu\in\Z_+^\ell$ with $\nu<\mu$,
suppose we have proved that $\phi(x^\a,x^{-\a}\ptl^\nu)=0$.
\vms{-4pt}Then
$$
\phi(x^\a,x^{-\a}\ptl^\mu)=
(\tau_\ell(\mu_\ell+1))^{-1}
\phi(x^\a,[x^{-\a-\tau}\ptl^{\mu+1_{[\ell]}},x^\tau])=0,
\vms{-4pt}\eqno(3.42)$$
where the first equality follows from the assumption that
$\phi(x^\a,x^{-\a}\ptl^\nu)=0$ for $\nu<\mu$,
and
the last equality from (1.2), (3.40), (3.41) and the assumption.
Finally,
$$
\matrix{
\phi(x^\a\ptl^\mu,x^{-\a}\ptl^\nu)
\!\!\!\!&=(\tau_\ell(\nu_\ell+1))^{-1}
\phi(x^\a\ptl^\mu,[x^{-\a-\tau}\ptl^{\nu+1_{[\ell]}},x^\tau])
\vs{4pt}\hfill\cr&
=(\tau_\ell(\nu_\ell+1))^{-1}
\phi(x^{-\a-\tau}\ptl^{\nu+1_{[\ell]}},[x^\a\ptl^\mu,x^\tau])
\vs{4pt}\hfill\cr&
=0,
\hfill\cr}
\eqno(3.43)$$
where the first equality follows from induction on $\mu+\nu$, and
the last equality follows from induction on $\mu$. Thus $\phi=0$.
\qed\par
We would like to
conclude this paper by defining the central extension $\wh{\WW}(\G)$
of the Lie algebra
$\WW(0,0,1,0,\G)$. Let $\phi_0$ be the 2-cocycle defined by (3.5) with
$\g=0$. Then $\wh{\WW}(\G)$ is the Lie algebra spanned by
$\{L_{\a,i,\mu},c
\,|\,(\a,i,\mu)\in\G\times\Z\times\Z_+\}$, with the following relations
\vms{-6pt}(cf.~(3.14), (3.5))
$$
\matrix{
[L_{\a,i,\mu},\!\!\!\!\!&L_{\b,j,\nu}]
\vs{2pt}\hfill\cr&\dis=
\dis\sum_{s=0}^{\mu+\nu}\sum_{r=0}^s
((^{\dis\mu}_{\ssc\,\dis s})(^{\dis s}_{\dis r})
[j]_r\b^{s-r}
-(^{\dis\nu}_{\ssc\,\dis s})(^{\dis s}_{\dis r})
[i]_r\a^{s-r})L_{\a+\b,i+j-r,\mu+\nu-s}
\vs{2pt}\hfill\cr&\dis\ \ \ \
+\d_{\a+\b,0}(-1)^\mu\mu!\nu!
\sum_{s=0}^{\mu+\nu+1}
(^{\ssc\,\dis i}_{\dis s})\d_{i+j,s-1}{\a^{\mu+\nu+1-s}\over(\mu+\nu+1-s)!}c,
\hfill\cr}
\vms{-6pt}\eqno(3.44)$$
and $[c,L_{\a,i,\mu}]=0$, for $(\a,i,\mu),(\b,j,\nu)\in\G\times\Z\times\Z_+$.
\par\
\par
\cl{                     ACKNOWLEDGEMENT
}
The author was supported by a fund from Education Department of China.
\par\
\par
\cl{                      REFERENCES
}
\par
\ni\hi3.4ex\ha1
[1] Berman S. On the low dimensional cohomology of some infinite
dimensional simple Lie algebras. Pac.~J.~Math. {\bf1979}, {\it83}, 27-36.
\par\ni\hi3.4ex\ha1
[2] Frenkel I.~B. Two constructions of affine Lie algebra representations
and boson-fermion correspondence in quantum field theory.
J.~Funct.~Anal. {\bf1981}, {\it44}, 259-327.
\par\ni\hi3.4ex\ha1
[3] Gelfand I.~M.; Fuks D.~B. The cohomology of the Lie algebra of
the vector fields in a circle. Funct.~Anal.~Appl.{1968}, {\it\bf2},
342-343.
\par\ni\hi3.4ex\ha1
[4] Kac V.~G. Lectures on infinite wedge representation and soliton
equations. Shuxue Jinzhan (Chinese Adv.~in Math.) {\bf1987}, {\it16},
343-376.
\par\ni\hi3.4ex\ha1
[5] Kac V.~G. {\it Infinite Dimensional Lie Algebras}, 3rd ed.;
Combridge Univ. Press, 1990.
\par\ni\hi3.4ex\ha1
[6] Kac V.~G.;  Kazhdan D.~A.; Lepowsky J.; Wilson R.~L.
Realization of the basic representation on the Euclidean Lie algebras.
adv.~in Math. {\bf1981}, {\it42}, 83-112.
\par\ni\hi3.4ex\ha1
[7] Kac V.~G.; Peterson D.~H. Spin and Wedge
representations of infinite dimensional Lie algebras and groups.
Proc.~Nat.~Acad.~Sci.~U.~S.~A. {\bf1981}, {\it78}, 3308-3312.
\par\ni\hi3.4ex\ha1
[8] Lepowsky J.; Wilson R.~L.
Constructions of the affine Lie algebra $A_1^{(1)}$.
Comm.~Math. Phys. {\bf1987}, {\it62}, 43-53.
\par\ni\hi3.4ex\ha1
[9] Li W. 2-Cocycles on the algebra of differential operators.
J.~Alg. {\bf1989}, {\it122}, 64-80.
\par\ni\hi4ex\ha1
[10] Li W.; Wilson R.~L. Central extensions of some Lie algebras.
Proc.~Amer. Math.~Soc. {\bf1998}, {\it126}, 2569-2577.
\par\ni\hi4ex\ha1
[11] Passman D.~P. Simple Lie algebras of Witt type. J.~Alg.
{\bf 1998}, {\it206}, 682-692.
\par\ni\hi4ex\ha1
[12] Patera J.; Zassenhaus H. The higher rank Virasoro algebras.
Comm. Math. Phys. {\bf1991}, {\it 136}, 1-14.
\par\ni\hi4ex\ha1
[13] Su, Y. On the low dimensional cohomology of Kac-Moody algebras
with coefficients in the complex field.
Shuxue Jinzhan (Chinese Adv.~in Math.) {\bf1989}, {\it18}, 346-351.
\par\ni\hi4ex\ha1
[14] Su, Y. 2-Cocycles on the Lie algebras of all differential operators
of several indeterminates.
(Chinese) Northeastern Math.~J. {\bf1990}, {\it6}, 365-368.
\par\ni\hi4ex\ha1
[15] Su, Y. A classification of indecomposable $sl_2(\C)$-modules
 and a conjecture of Kac on Irreducible modules over the Virasoro
 algebra. J. Alg. {\bf1993}, {\it 161}, 33-46.
\par\ni\hi4ex\ha1
[16] Su, Y. Harish-Chandra modules of the intermediate series over
the high rank Virasoro algebras and high rank super-Virasoro
algebras. J. Math. Phys. {\bf1994}, {\it 35}, 2013-2023.
\par\ni\hi4ex\ha1
[17] Su, Y.; Zhao K.
Generalized Virasoro and super-Virasoro algebras and modules of the
intermediate series. To appear.
\par\ni\hi4ex\ha1
[18] Su, Y.; Zhao K.
Second cohomology group of generalized Witt type Lie algebras and
certain representations.
To appear.
\par\ni\hi4ex\ha1
[19] Su, Y.; Zhao K. Simple algebras of Weyl types.
Science in China A. In press.
\par\ni\hi4ex\ha1
[20]  Su, Y.; Zhao K. Isomorphism classes and automorphism groups of
algebras of Weyl type, Science in China A. In press.
\par\ni\hi4ex\ha1
[21] Xu X.
New generalized simple Lie algebras of Cartan type over a field with
characteristic 0. J. Alg. {\bf2000}, {\it224}, 23-58.
\end{document}